# Probability measures, Lévy measures and analyticity in time

OLE E. BARNDORFF-NIELSEN[1] and FRIEDRICH HUBALEK[2]

[1]*Department of Mathematical Sciences, University of Aarhus, Ny Munkegade, DK-8000 Aarhus C, Denmark*
[2]*Institute of Mathematical Methods in Economics, Vienna University of Technology, Wiedner Hauptstrasse 8/105-1, A-1040 Vienna, Austria. E-mail: fhubalek@fam.tuwien.ac.at*

We investigate the relation of the semigroup probability density of an infinite activity Lévy process to the corresponding Lévy density. For subordinators, we provide three methods to compute the former from the latter. The first method is based on approximating compound Poisson distributions, the second method uses convolution integrals of the upper tail integral of the Lévy measure and the third method uses the analytic continuation of the Lévy density to a complex cone and contour integration. As a by-product, we investigate the smoothness of the semigroup density in time. Several concrete examples illustrate the three methods and our results.

## 1. Introduction

For the infinitely divisible laws, there are a number of intriguing and useful relations and points of similarity between the probability measures or probability densities of the laws on the one hand and their associated Lévy measures or Lévy densities on the other.

In particular, if $U$ is the Lévy measure of an infinitely divisible law on $\mathbb{R}^d$ with associated Lévy process $\{X_t\}_{t \geq 0}$ and if $P(\mathrm{d}x; t)$ denotes the law of $X_t$, then (see Sato (1999), Corollary 8.9)

$$\lim_{t \to 0} t^{-1} \int_{\mathbb{R}^d} f(x) P(\mathrm{d}x; t) = \int_{\mathbb{R}^d} f(x) U(\mathrm{d}x) \tag{1}$$

for any function $f$ in the space $C_{\#}$ of bounded continuous functions on $\mathbb{R}^d$ vanishing in a neighborhood of 0. It may also be noted that Burnaev (2006) gives two formulas for computing the Lévy measure from the corresponding cumulant function.

The present paper considers the opposite problem, that of calculating $P(\mathrm{d}x; t)$ from $U$.

---

In many concrete examples, there is a power series representation

$$P(\mathrm{d}x; t) = \sum_{n=1}^{\infty} \frac{t^n}{n!} U_n(\mathrm{d}x), \tag{2}$$

where the $U_n(\mathrm{d}x)$ are, in general, signed measures and $U_1$ is equal to $U$, the Lévy measure of $P(\mathrm{d}x; t)$. When $d = 1$ and both $P(\mathrm{d}x; t)$ and $U(\mathrm{d}x)$ are concentrated on the positive half-line $\mathbb{R}_{>0} = (0, \infty)$, it is convenient to give (2) the form

$$P^+(x; t) = \sum_{n=1}^{\infty} \frac{t^n}{n!} U_n^+(x), \tag{3}$$

with $P^+(x; t) = P((x, +\infty); t)$ and $U_n^+(x) = U_n((x, +\infty))$ denoting the upper tail integrals.

More particularly, when $P(\mathrm{d}x, t)$ (resp. $U(\mathrm{d}x)$) is absolutely continuous with density $p(x; t)$ resp. $u(x)$ – the setting on which we focus in this paper – then the 'density version' of (3) is

$$p(x; t) = \sum_{n=1}^{\infty} \frac{t^n}{n!} u_n(x). \tag{4}$$

In (4), necessarily, $u_1 = u$ and the question is how the further coefficients $u_n$ may be *calculated* from $u$, possibly also using properties of the cumulant function of $X_1$ (which, of course, is essentially determined by $u$). As part of the problem, we discuss conditions ensuring that an expansion of the form (4) *exists*. This issue is of some independent interest.

Note that in the case of a finite Lévy measure, the process $X$ is a compound Poisson process and validity of formulae (4) and (2), with straightforward modifications for the atom at zero, is easily established; see Barndorff-Nielsen and Hubalek (2006), Section 4.1. We shall not consider this case here any further.

Except for the discussion in Section 2 and some remarks in Section 4, we only consider the case where the process $X$ is a subordinator with infinite Lévy measure and without linear drift. In other words, $X$ is an infinite activity pure jump subordinator.

We shall discuss three methods for determining the coefficients $u_n(x)$. The first involves, as the final step, a limiting operation. We have $u_n(x) = \lim_{\varepsilon \to 0} u_{n\varepsilon}(x)$, where

$$u_{n\varepsilon}(x) = \sum_{k=1}^{n} (-1)^{n-k} \binom{n}{k} c(\varepsilon)^{n-k} u_\varepsilon^{*k}(x). \tag{5}$$

Here, $u_\varepsilon(x)$ is an approximation of the Lévy density $u(x)$ that corresponds to a compound Poisson process with intensity $c(\varepsilon)$ and $*k$ indicates $k$-fold convolution. The second method uses derivatives of convolutions, namely,

$$u_n(x) = (-1)^n \frac{\mathrm{d}^n}{\mathrm{d}x^n} ((U^+)^{*n})(x). \tag{6}$$

In words: we can obtain $u_n(x)$ as the $n$th derivative of the $n$th convolution power of the upper tail integral of the Lévy density $u(x)$.

The third method uses the complex contour integral

$$u_n(x) = \frac{1}{2\pi i} \int_{\mathcal{C}} \kappa(\theta)^n \mathrm{e}^{\theta x}\, \mathrm{d}\theta, \tag{7}$$

where $\kappa(\theta)$ is the analytic continuation of the cumulant function to a complex cone containing $\Re\theta \geq 0$ and the contour is, roughly speaking, along the boundary of the cone. We will see that such an analytic continuation can be derived from an analytic continuation of the Lévy density $u(x)$ to a complex cone containing the positive real axis.

We proceed to mention various works containing results that are related to those of the present paper. The comprehensive monograph by Sato (1999) contains many instances of the interesting relations between the probability distributions and the Lévy measures of infinitely divisible laws (cf. also Embrechts *et al.* (1979), Embrechts and Goldie (1981) and Sato and Steutel (1998)). Some examples are the relation between unimodality properties of the two types of densities (Sato (1999), Section 52) and the behavior under exponential tilting (or Esscher transformation). See also Léandre (1987), Ishikawa (1994) and Picard (1997) who, partly in the wider setting of pure jump processes, study cases where the transition density exists and behaves as a power of $t$ for $t \to 0$. Continuity of $P^+(x; t)$ at $t = 0$ is characterized in Doney (2004). In Rüschendorf and Woerner (2002) (also see Woerner (2001)), the authors have established the validity of expansions for the probability density or distribution function of $X_t$ that are related to, but essentially different from, (3) and (4).

The paper is organized as follows. Section 2 consists of a number of initial remarks on the problem at hand. Section 3 contains our main mathematical results. Illustrative examples will be given in Section 4. Technical auxiliary material used in the proofs appears in the Appendix.

The results in the present paper build partly on our previous, unpublished preliminary work (Barndorff-Nielsen (2000) and Hubalek (2002)).

## 2. Initial remarks

### 2.1. A first, motivating example: the positive $\alpha$-stable distribution

We consider the positive $\alpha$-stable distribution with Lévy density

$$u(x) = -\frac{x^{-1-\alpha}}{\Gamma(-\alpha)}, \tag{8}$$

where $0 < \alpha < 1$. Note that we interpret $\Gamma(s)^{-1}$ as an entire function with zeros at the non-positive integers. Alternatively, we could use the functional equations $\Gamma(s + 1) =$ $s\Gamma(s)$ and $\Gamma(s)\Gamma(1-s) = \pi\csc(\pi s)$ to rewrite expressions in a more familiar (and more lengthy) form.

In general, there is no closed form expression for $p(x;t)$ in terms of elementary functions, but it is well known (see, e.g., Feller (1971), XVII.7 (6.8)) that

$$p(x;t) = \sum_{n\geq 1}\frac{t^n}{n!}, \qquad u_n(x) = \frac{(-1)^n}{\Gamma(-n\alpha)}x^{-1-n\alpha}. \tag{9}$$

Let us first illustrate the calculation of $u_2(x)$ by our first method, that is, as limit of $u_{2\varepsilon}(x)$ for $\varepsilon \to 0$, where we use the approximation $u_\varepsilon(x) = I_{[\varepsilon,\infty)}(x)u(x)$. We have $c(\varepsilon) = \varepsilon^{-\alpha}/\Gamma(1-\alpha)$. For $x > 2\varepsilon$, we obtain, by symmetry and partial integration,

$$u_\varepsilon^{*2}(x) = \frac{2}{\alpha\Gamma(-\alpha)^2}\Bigg[\varepsilon^{-\alpha}(x-\varepsilon)^{-1-\alpha}$$
$$- \left(\frac{x}{2}\right)^{-1-2\alpha} + (1+\alpha)\int_\varepsilon^{x/2} y^{-\alpha}(x-y)^{-2-\alpha}\,\mathrm{d}y\Bigg]. \tag{10}$$

As $u_{2\varepsilon}(x) = u_\varepsilon^{*2}(x) - 2c(\varepsilon)u_\varepsilon(x)$, we obtain, in the limit,

$$u_2(x) = \frac{2}{\alpha\Gamma(-\alpha)^2}\Bigg[-\left(\frac{x}{2}\right)^{-1-2\alpha} + (1+\alpha)\int_0^{x/2} y^{-\alpha}(x-y)^{-2-\alpha}\,\mathrm{d}y\Bigg]. \tag{11}$$

The integral on the right-hand side can be expressed in terms of the incomplete beta function and, in this particular case, reduced to integrals for the complete beta function by elementary substitutions. This finally yields agreement with (9) for $n=2$.

In principle, though less explicit and more cumbersome, the method can be used for $n>2$. Instead, let us illustrate our second method, the calculation of $u_n(x)$ according to formula (6). The tail integral is $U^+(x) = x^{-\alpha}/\Gamma(1-\alpha)$ and, by induction, or, quicker, by looking at the Laplace transforms, we see that $(U^+)^{*n}(x) = x^{n-n\alpha}/\Gamma(n-n\alpha)$. Differentiating this equation $n$ times and applying the functional equation of the Gamma function to simplify the expression, we obtain (9).

Finally, let us illustrate our third method, the calculation of $u_n(x)$ according to formula (7). The Laplace cumulant function is $\kappa(\theta) = -\theta^\alpha$ and we get

$$u_n(x) = \frac{(-1)^n}{2\pi\mathrm{i}}\int_C \theta^{n\alpha}\mathrm{e}^{\theta x}\,\mathrm{d}\theta. \tag{12}$$

To see that this, in fact, gives (9), we have to substitute $\theta \mapsto \theta/x$ and recognize the resulting integral as a variant of the Hankel contour integral for $\Gamma(-n\alpha)^{-1}$.

## 2.2. A simple, general result

Let $P(x;t)$ and $u(x)$ be, respectively, the cumulative distribution function and the Lévy density (assumed to exist) of an infinite activity Lévy process on $\mathbb{R}_{>0}$. Let the $u_\varepsilon(x)$ be integrable Lévy densities that we think of as approximations of $u(x)$ and let us define $c(\varepsilon)$ and $u_{n\varepsilon}(x)$ as in (5) above, setting

$$U_{0\varepsilon}(x) = 1, \qquad U_{n\varepsilon}(x) = -\int_x^\infty u_{n\varepsilon}(y)\,\mathrm{d}y \qquad (n \geq 1). \tag{13}$$

**Theorem 1.** *Suppose*

$$\lim_{\varepsilon \to 0} \int (1 \wedge x)|u_\varepsilon(x) - u(x)|\,\mathrm{d}x = 0. \tag{14}$$

*Then*

$$P(x;t) = \lim_{\varepsilon \to 0} \sum_{n \geq 0} U_{n\varepsilon}(x)\frac{t^n}{n!}, \tag{15}$$

*pointwise for each* $x \in \mathbb{R}_{>0}$ *and* $t > 0$.

This follows easily from Lévy's continuity theorem and the observation that the distribution of an infinite activity process is continuous. Details can be found in Barndorff-Nielsen (2000) and Barndorff-Nielsen and Hubalek (2006). Note that the approximation $u_\varepsilon(x)I_{[\varepsilon,\infty)}$ satisfies the assumptions of Theorem 1.

## 2.3. Formulae for $u_{n\varepsilon}$

The first few functions $u_{n\varepsilon}$ are

$$\begin{aligned}
u_{1\varepsilon}(x) &= u_\varepsilon(x), \\
u_{2\varepsilon}(x) &= u_\varepsilon^{*2}(x) - 2c(\varepsilon)u_\varepsilon(x), \\
u_{3\varepsilon}(x) &= u_\varepsilon^{*3}(x) - 3c(\varepsilon)u_\varepsilon^{*2}(x) + 3c(\varepsilon)^2 u_\varepsilon(x), \\
u_{4\varepsilon}(x) &= u_\varepsilon^{*4} - 4c(\varepsilon)u_\varepsilon^{*3} + 6c(\varepsilon)^2 u_\varepsilon^{*2}(x) - 4c(\varepsilon)^3 u_\varepsilon(x).
\end{aligned} \tag{16}$$

Further, it follows from the well-known inverse relations for binomial sums (see, e.g., Comtet (1970), III.6.a e) that the formula defining $u_{n\varepsilon}(x)$ in (5) implies conversely that

$$u_\varepsilon^{*n}(x) = \sum_{k=1}^n \binom{n}{k} c(\varepsilon)^{n-k} u_{k\varepsilon}(x). \tag{17}$$

This can be used to compute $u_{n\varepsilon}(x)$ inductively by

$$u_{n\varepsilon}(x) = u_\varepsilon^{*n}(x) - \sum_{k=1}^{n-1} \binom{n}{k} c(\varepsilon)^{n-k} u_{k\varepsilon}(x). \tag{18}$$

In particular, we have

$$
\begin{aligned}
u_{1\varepsilon}(x) &= u_{\varepsilon}(x), \\
u_{2\varepsilon}(x) &= u_{\varepsilon}^{*2}(x) - 2c(\varepsilon)u_{1\varepsilon}(x), \\
u_{3\varepsilon}(x) &= u_{\varepsilon}^{*3}(x) - 3c(\varepsilon)u_{2\varepsilon}(x) - 3c(\varepsilon)^2 u_{\varepsilon}(x), \\
u_{4\varepsilon}(x) &= u_{\varepsilon}^{*4} - 4c(\varepsilon)u_{3\varepsilon} - 6c(\varepsilon)^2 u_{2\varepsilon}(x) - 4c(\varepsilon)^3 u_{\varepsilon}(x).
\end{aligned}
\tag{19}
$$

## 2.4. Cancellation of singularities

Convergence of $u_{n\varepsilon}(x)$ to a function $u_n(x)$ implies a subtle cancellation of singularities. Equation (5) contains an alternating sum of terms that diverge. Both $c(\varepsilon)$ and the convolution powers of $u_{\varepsilon}(x)$ tend to $+\infty$ as $\varepsilon \to 0$; in particular, see the first few instances of that formula listed in (16) above.

To gain an understanding of how this cancellation occurs, note that for $n \geq 0$, we have

$$
u_{n+1\varepsilon}(x) = (n+1)x^{-1}\left\{\int_0^x u_{n\varepsilon}(x-y)\bar{u}_{\varepsilon}(y)\,\mathrm{d}y + (-1)^n c(\varepsilon)^n \bar{u}_{\varepsilon}(x)\right\},
\tag{20}
$$

with $\bar{u}_{\varepsilon}(x) = xu_{\varepsilon}(x)$, as can be established by induction. Consider the case $n = 2$, and let $U_{\varepsilon}^+(x) = \int_x^{\infty} u_{\varepsilon}(y)\,\mathrm{d}y$ and $U^+(x) = \int_x^{\infty} u(y)\,\mathrm{d}y$. Using (20) and noting that $c(\varepsilon) = \int_0^x u_{\varepsilon}(y)\,\mathrm{d}y + U_{\varepsilon}^+(x)$, we may rewrite $u_{2\varepsilon}(x)$ as

$$
u_{2\varepsilon}(x) = 2x^{-1}\left\{\int_0^x u_{\varepsilon}(y)\{\bar{u}_{\varepsilon}(x-y) - \bar{u}_{\varepsilon}(x)\}\,\mathrm{d}y - \bar{u}_{\varepsilon}(x)U_{\varepsilon}^+(x)\right\}.
\tag{21}
$$

Hence, by letting $\varepsilon \to 0$ and invoking condition (14), we obtain the following.

**Proposition 2.** *Suppose the Lévy density $u(x)$ is differentiable and let $\bar{u}(x) = xu(x)$. Then*

$$
u_2(x) = 2x^{-1}\left\{\int_0^x u(y)\{\bar{u}(x-y) - \bar{u}(x)\}\,\mathrm{d}y - \bar{u}(x)U^+(x)\right\},
\tag{22}
$$

*with the integral existing and being finite.*

Formula (21), first given in Barndorff-Nielsen (2000), has been generalized by Woerner (2001) to

$$
\frac{1}{n+1}\bar{u}_{n+1\varepsilon}(x) = \frac{1}{n}\left[\int_0^x u_{\varepsilon}(y)\{\bar{u}_{n\varepsilon}(x-y) - \bar{u}_{n\varepsilon}(x)\}\,\mathrm{d}y - \bar{u}_{n\varepsilon}(x)U^+(x)\right],
\tag{23}
$$

where $\bar{u}_{n\varepsilon}(x) = xu_{n\varepsilon}(x)$. Typically, though, this cannot be used to pass to the limit $\varepsilon \to 0$ for $n > 1$.

## 2.5. Miscellaneous further points

- The power series representation (4) with the coefficients $u_n(x)$ given as limit of $u_{n\varepsilon}(x)$ was originally derived by a heuristic argument that may be found in Barndorff-Nielsen and Shephard (2008), Section 3.3.3; see also Barndorff-Nielsen and Hubalek (2006), Section 2.1. In Section 3 below, we give a rigorous derivation for the case of one-dimensional subordinators, under fairly strong assumptions. In subsequent work, we hope to establish proofs for Lévy processes in $\mathbb{R}^d$.

- Heavy tails do not matter much in the problems studied in the present context as we can deal with them using the Esscher transform. Rather, it is the behavior of small jumps that can cause difficulties, as reflected in the assumptions of the theorems given below. For details, see Barndorff-Nielsen and Hubalek (2006), Section 2.4.

- Suppose $p(x; 0) = 0$ and

$$\lim_{t \to 0} t^{-1} p(x; t) = u(x), \tag{24}$$

and let $u_\varepsilon(x) = \varepsilon^{-1} p(x; \varepsilon)$, in which case $c(\varepsilon) = \varepsilon^{-1}$. By using the semigroup property $p^{*k}(x, \varepsilon) = p(x, k\varepsilon)$, formula (5) takes the form

$$u_{n\varepsilon}(x) = \varepsilon^{-n} \sum_{k=1}^{n} (-1)^{n-k} \binom{n}{k} p(x; k\varepsilon). \tag{25}$$

Now, the right-hand side of (25) is, in fact, an $n$th order difference quotient of $p(x; t)$, so, provided $p(x; t)$ is $n$ times differentiable from the right at $t = 0$, we have

$$\lim_{\varepsilon \to 0} u_{n\varepsilon}(x) = \frac{\partial^n}{\partial t^n} p(x; 0). \tag{26}$$

Thinking of $u_\varepsilon(x)$ as an approximation of $\varepsilon^{-1} p(x; \varepsilon)$, as well as an approximation of $u(x)$, this is then a further indication that, in considerable generality, univariate as well as multivariate, $p(x; t)$ may be calculated via (5) in the manner discussed above. Of course, in practice, choosing $u_\varepsilon(x) = \varepsilon^{-1} p(x; \varepsilon)$ is not an option since the point is to determine $p(x; t)$ in terms of the Lévy density $u(x)$.

- As a referee pointed out to us, our work has a strong connection to semigroup theory, in particular to the so-called *exponential formulae*. Let $T$ be the transition operator defined by

$$Tf(x) = \int_0^\infty f(x + y) p(y; t) \, \mathrm{d}y. \tag{27}$$

For example, Theorem 1, which holds true for infinitely divisible distributions on $\mathbb{R}^d$, can be seen as a variant of Hille and Phillips (1957) $(E_2)$, page 354. Yet, our main interest is in interchanging the delicate limit $\varepsilon \to 0$ and the infinite summation, and we have not been able to use results from semigroup theory for that purpose.

- Condition (14) is satisfied, in particular, if there exists an integrable function $v$ on $\mathbb{R}_{>0}$ such that $(x \wedge 1)u_\varepsilon(x) \leq v(x)$ for all $x \in \mathbb{R}_{>0}$ and all $\varepsilon$. Some candidates for $u_\varepsilon$ are

$$u_\varepsilon(x) = \mathbf{1}_{[\varepsilon,\infty)}(x)u(x) \quad \text{or} \quad u_\varepsilon(x) = u(x)\mathrm{e}^{-x/\varepsilon}. \tag{28}$$

# 3. Main results

In this section, we analyze the following issues.

- Does $u_{n\varepsilon}(x)$ converge as $\varepsilon \to 0$? If so, can we find a more direct method to compute the limit $u_n(x)$ from $u(x)$ and so avoid the difficult cancellations in $u_{n\varepsilon}(x)$ as $\varepsilon \to 0$?
- If we have convergence, is $p(x;t)$ in fact $n$-times differentiable (from the right) at $t = 0$ and, if so, is $u_n(x)$ the $n$th derivative?
- If the answer to the previous question is yes for all $n \geq 1$, do we have a convergent Taylor expansion of $p(x;t)$ at $t = 0$? Is $p(x;t)$ in fact an entire function in $t \in \mathbb{C}$?

In the proofs, we refer to several technical estimates that are provided in the form of lemmas in the Appendix. Recall that we are assuming that the process $X$ is an infinite activity subordinator.

## 3.1. Pointwise convergence of the coefficient functions

In this subsection, we investigate the limiting behavior of $u_{n\varepsilon}(x)$ as $\varepsilon \to 0$ for the particular choice $u_\varepsilon(x) = \mathrm{e}^{-\varepsilon/x}u(x)$. This approximation is simple, it is always feasible, in the sense that it implies

$$\lim_{\varepsilon \to 0} u_\varepsilon(x) = u(x), \qquad \lim_{\varepsilon \to 0} \int_x^\infty u_\varepsilon(y)\,\mathrm{d}y = \int_x^\infty u(y)\,\mathrm{d}y \qquad \forall x > 0, \tag{29}$$

and $u_\varepsilon(x)$ will be smooth if $u(x)$ is smooth, a property exploited below. We provide conditions on $u(x)$ that imply the convergence of $u_{n\varepsilon}(x)$ and obtain the expression (6) for the limit $u_n(x)$.

**Theorem 3.** *Suppose $n \in \mathbb{N}$ and*

$$\int_0^\infty \mathrm{e}^{-rx}x^{k+1}|u^{(k)}(x)|\,\mathrm{d}x < \infty, \qquad k = 0, 1, \ldots, m, \tag{30}$$

*holds for some $m \geq n+2$ and $r > 0$.*

(i) *Let*

$$U^+(x) = \int_x^\infty u(y)\,\mathrm{d}y, \qquad x > 0. \tag{31}$$

*Then the $n$th convolution power of $U^+(x)$ is well defined for $x > 0$ and $(U^+)^{*n}(x)$ is $n$-times continuously differentiable.*

(ii) *If we set $u_\varepsilon(x) = \mathrm{e}^{-\varepsilon/x} u(x)$ for $\varepsilon > 0$ and $x > 0$ and define*

$$u_{n\varepsilon}(x) = \sum_{k=1}^{n} (-1)^{n-k} \binom{n}{k} u_\varepsilon^{*k}(x) c(\varepsilon)^{n-k}, \tag{32}$$

*where $c(\varepsilon) = \int_0^\infty u_\varepsilon(x)\,\mathrm{d}x$, then $\lim_{\varepsilon \to 0} u_{n\varepsilon}(x) = u_n(x)$ for $x > 0$, where*

$$u_n(x) = (-1)^n \frac{\mathrm{d}^n}{\mathrm{d}x^n} (U^+)^{*n}(x). \tag{33}$$

**Proof.** In Lemma A.9, we show that $(U^+)^{*n}(x)$ exists, and in Lemma A.10, that it is $n$-times differentiable. Smoothness of convolutions is less obvious than it at first seems (cf. Doetsch (1950), Section 2.14, page 104ff and Uludağ (1998)).

Let $\lambda_{n\varepsilon}(\theta)$ denote the Laplace transform of $u_{n\varepsilon}(x)$. It is given by $\lambda_{n\varepsilon}(\theta) = \kappa_\varepsilon(\theta)^n - (-1)^n c(\varepsilon)^n$. Let $\lambda_n(\theta) = \kappa(\theta)^n$. Note that $\lambda_n(\theta)$ is *not* the Laplace transform of $u_n(x)$; the Laplace transform of $u_n(x)$ does not exist. But, using the estimates from Lemma A.5, namely

$$|\lambda_{n\varepsilon}^{(m)}(\theta)| \leq E_{mn}/|\theta|^{m-n}, \qquad |\lambda_n^{(m)}(\theta)| \leq E_{mn}/|\theta|^{m-n} \tag{34}$$

for some constants $E_{mn}$, we see that both $x^m u_{n\varepsilon}(x)$ and $x^m u_n(x)$ have integrable Laplace transforms, namely $\lambda_{n\varepsilon}^{(m)}(\theta)$ and $\lambda_n^{(m)}(\theta)$. According to Lemma A.10, we can write the inversion integrals

$$u_{n\varepsilon}(x) = \frac{1}{2\pi\mathrm{i}x^m} \int_{r-\mathrm{i}\infty}^{r+\mathrm{i}\infty} \lambda_{n\varepsilon}^{(m)}(\theta) \mathrm{e}^{\theta x}\,\mathrm{d}\theta, \qquad u_n(x) = \frac{1}{2\pi\mathrm{i}x^m} \int_{r-\mathrm{i}\infty}^{r+\mathrm{i}\infty} \lambda_n^{(m)}(\theta) \mathrm{e}^{\theta x}\,\mathrm{d}\theta, \tag{35}$$

and obtain

$$|u_{n\varepsilon}(x) - u_n(x)| \leq \frac{\mathrm{e}^{rx}}{2\pi x^m} \int_{-\infty}^{+\infty} |\lambda_{n\varepsilon}^{(m)}(r+\mathrm{i}y) - \lambda_n^{(m)}(r+\mathrm{i}y)|\,\mathrm{d}y. \tag{36}$$

We have $0 \leq u_\varepsilon(x) \leq u(x)$, for $x > 0, \varepsilon > 0$ and $u_\varepsilon(x) \to u(x)$ for $\varepsilon \to 0$. Thus, looking at

$$\kappa_\varepsilon(\theta) = \int_0^\infty (\mathrm{e}^{-\theta x} - 1) u_\varepsilon(x)\,\mathrm{d}x, \qquad \kappa(\theta) = \int_0^\infty (\mathrm{e}^{-\theta x} - 1) u(x)\,\mathrm{d}x \tag{37}$$

and

$$\kappa_\varepsilon^{(k)}(\theta) = (-1)^k \int_0^\infty \mathrm{e}^{-\theta x} x^k u_\varepsilon(x)\,\mathrm{d}x, \qquad \kappa^{(k)}(\theta) = (-1)^k \int_0^\infty \mathrm{e}^{-\theta x} x^k u(x)\,\mathrm{d}x \tag{38}$$

for $k = 1, \ldots, n$, we see, by dominated (or monotone) convergence, that $\lim_{\varepsilon \to 0} \kappa_\varepsilon(\theta) = \kappa(\theta)$, $\lim_{\varepsilon \to 0} \kappa_\varepsilon^{(k)}(\theta) = \kappa^{(k)}(\theta)$. The functions $\lambda_{n\varepsilon}^{(m)}(\theta)$ and $\lambda_n^{(m)}(\theta)$ are polynomials in $\kappa_\varepsilon(\theta)$, $\kappa_\varepsilon'(\theta), \ldots, \kappa_\varepsilon^{(m)}(\theta)$, respectively in $\kappa(\theta)$, $\kappa'(\theta), \ldots, \kappa^{(m)}(\theta)$, thus $\lim_{\varepsilon \to 0} \lambda_{n\varepsilon}^{(m)}(\theta) = \lambda_n^{(m)}(\theta)$. Moreover, they are dominated by the integrable function $E_{mn}/|\theta|^{m-n}$ and by dominated convergence in (36), we have the desired result. $\qquad \square$

An interesting class of infinitely divisible distributions on $\mathbb{R}_{>0}$ is the family of generalized gamma convolutions. This class is characterized by having absolutely continuous Lévy measures with densities $u(x)$ such that $\bar{u}(x) = xu(x)$ are completely monotone functions (Bondesson (1992), Theorem 3.1.1).

**Theorem 4.** *If $u(x)$ is a Lévy density such that $\bar{u}(x) = xu(x)$ is completely monotone, then the integrability assumptions (30) in Theorem 3 hold for all $n \in \mathbb{N}$ with arbitrary $r > 0$.*

**Proof.** By the Bernstein–Widder representation for completely monotone functions, we know $\bar{u}(x)$ is holomorphic in $\Re x > 0$, thus the Taylor series at $x > 0$ has radius of convergence $x$ and

$$\bar{u}(x/2) = \bar{u}(x - x/2) = \sum_{n=0}^{\infty} (-1)^n \bar{u}^{(n)}(x) \frac{x^n}{2^n n!}. \tag{39}$$

In our setting, $\int_0^\infty \mathrm{e}^{-\theta x} \bar{u}(x/2) \, \mathrm{d}x < \infty$ for $\theta > 0$ and as $(-1)^n \bar{u}^{(n)}(x) \geq 0$, we can integrate the series term by term. As $\bar{u}^{(n)}(x) = xu^{(n)}(x) + nu^{(n-1)}(x)$ for $n \geq 1$ and we know that $\int \mathrm{e}^{-rx} xu(x) \, \mathrm{d}x < \infty$, we inductively obtain the result. $\qquad\square$

**Remark 5.** An example where Theorem 3 applies but the corresponding distribution is not a generalized gamma convolution is given by $u(x) = x^{-3/2} \mathrm{e}^{\sin(x)}$. An example where the integrability conditions (30) do not hold for $n \geq 1$ and any $r > 0$ is given by $u(x) = x^{-3/2} \sin(x^{-3})^2$. We do not know whether the conclusion of the theorem is nevertheless true in this case.

## 3.2. Differentiability in time

For the proof of differentiability properties of the probability densities $p(x; t)$ with respect to $t \geq 0$, we need slightly different integrability properties of the cumulant function $\kappa(\theta)$ and its derivatives. Sufficient conditions to guarantee those from assumptions on the Lévy density $u(x)$ are conveniently formulated in terms of the *integral modulus of continuity*. We recall, for example, from Kawata (1972), Theorem 2.7.4, that the integral modulus of continuity $\omega^{(1)}(\delta; f)$ for an integrable function $f(x)$ and a real number $\delta > 0$ is defined by

$$\omega^{(1)}(\delta; f) = \sup_{0 < |h| \leq \delta} \int_{-\infty}^{+\infty} |f(x+h) - f(x)| \, \mathrm{d}x. \tag{40}$$

We note that it is sufficient to consider $0 < h \leq \delta$ in (40). We use the integral modulus of continuity for functions $f(x)$ that are a priori defined for $x > 0$ with the understanding that $f(x) = 0$, if $x \leq 0$.

**Theorem 6.** *Suppose $m \in \mathbb{N}$, $n \in \mathbb{N}$, $\alpha \in (0,1)$, $r \in [0,\infty)$ and $u(x)$ is the Lévy density of an infinite activity subordinator. Suppose*

$$m > \frac{1 + n\alpha}{1 - \alpha}, \tag{41}$$

*$u(x)$ is $m$-times differentiable in $x > 0$, the functions*

$$v_\ell(x) = (-1)^\ell \mathrm{e}^{-rx} x^{\ell+1} u^{(\ell)}(x) \qquad (\ell = 0, \ldots, m) \tag{42}$$

*are integrable and their integral modulus of continuity satisfies*

$$\omega^{(1)}(\delta; v_\ell) = \mathcal{O}(\delta^{1-\alpha}) \qquad (\delta \to 0). \tag{43}$$

*Let $p(x;t)$ denote the probability densities corresponding to $u(x)$. Then $p(x;t)$ is, for all $x > 0$, $n$-times differentiable in $t \geq 0$; furthermore,*

$$u_k(x) = (-1)^k \frac{\partial^k}{\partial x^k} (U^+)^{*k}(x) \qquad (k = 1, \ldots, n), \tag{44}$$

*is well defined and*

$$\frac{\partial^k}{\partial t^k} p(x, 0) = u_k(x) \qquad (k = 1, \ldots, n). \tag{45}$$

**Proof.** Let

$$\lambda_n^{(m)}(\theta;t) = \frac{\partial^{m+n}}{\partial \theta^m \, \partial t^n} \mathrm{e}^{\kappa(\theta)t}, \qquad p_n(x;t) = \frac{\partial^n}{\partial t^n} p(x;t). \tag{46}$$

We will show the following statement inductively for $n' = 1, \ldots, n$. We have, for all $x > 0$, that $p(x;t)$ is $n'$-times differentiable in $t \geq 0$ and that

$$p_{n'}(x;t) = \frac{(-1)^m}{2\pi \mathrm{i} x^m} \int_{c-\mathrm{i}\infty}^{c+\mathrm{i}\infty} \lambda_{n'}^{(m)}(\theta;t) \mathrm{e}^{\theta x} \, \mathrm{d}\theta. \tag{47}$$

First, with $n' = 0$, we observe that $(-1)^m x^m p(x;t)$ has Laplace transform $\lambda_0^{(m)}(\theta;t)$. In Lemma A.8 below, we show that the assumptions on the integral modulus of continuity imply $\lambda_0^{(m)}(\theta;t) = \mathcal{O}(|\theta|^{-m(1-\alpha)})$ as $\Im(\theta) \to \pm\infty$. Since (41) implies $m > 1/(1-\alpha)$, we have that $\lambda_0^{(m)}(\theta;t)$ is integrable and we can apply the inversion formula. This is all that was to be shown for $n' = 0$. Suppose, now, we have shown the claim for some $n' - 1$ and want to show it for $n'$. We can write

$$\begin{aligned}
&h^{-1}(p_{n'-1}(x;t+h) - p_{n'-1}(x;t)) \\
&= \frac{(-1)^m}{2\pi \mathrm{i} x^m} \int_{c-\mathrm{i}\infty}^{c+\mathrm{i}\infty} h^{-1}(\lambda_{n'-1}^{(m)}(\theta;t+h) - \lambda_{n'-1}^{(m)}(\theta;t)) \mathrm{e}^{\theta x} \, \mathrm{d}\theta.
\end{aligned} \tag{48}$$

In Lemma A.8 below, we show that the assumptions on the integral modulus of continuity imply $h^{-1}(\lambda_{n'-1}^{(m)}(\theta; t+h) - \lambda_{n'-1}^{(m)}(\theta; t)) = \mathcal{O}(|\theta|^{\alpha(m+n')-m})$ as $\Im(\theta) \to \pm\infty$. Now, (41) implies that the integrand in (48) is dominated by an integrable function and we can apply dominated convergence as $h \to 0$. This shows that $p_{n'}(x; t)$ is differentiable with respect to $t$ and its derivative is given by (47). This finishes the induction. To complete the proof, we observe that $\lambda_n(\theta) = \kappa(\theta)^n = \frac{\partial^n}{\partial t^n}[e^{\kappa(\theta)t}]_{t=0} = \lambda_n^{(0)}(\theta; 0)$ and, in view of (46), (47) and (35), we indeed obtain $p_n(x; 0) = u_n(x)$.  □

**Remark 7.** If we consider $n = 1$, then Theorem 6 provides sufficient conditions for the (pointwise) validity of

$$u(x) = \lim_{t \to 0} t^{-1} p(x; t) \tag{49}$$

for $x > 0$. Our assumptions are quite different from those given by Woerner (2001), Corollary 2.3.

### 3.3. Power series representation in time

The purpose of this section is to show that, subject to some regularity conditions, the probability densities $p(x; t)$ are *analytic functions* in $t$ that can be represented by a power series of the form (4). To be able to do so, we assume that the Lévy density $u(x)$ is an analytic function satisfying some growth condition.

**Theorem 8.** *Suppose*

$$a > 0, \qquad 0 < \alpha < 1, \qquad \beta > -1, \qquad \gamma > 0, \qquad 0 < \psi < \frac{\pi}{2} \tag{50}$$

*and the Lévy density $u(z)$ is an analytic function in a domain containing*

$$W = \{z \in \mathbb{C} : z \neq 0, |\arg(z)| \leq \psi\}. \tag{51}$$

*Assume, moreover, that*

$$u(z) = az^{-1-\alpha} + \mathcal{O}(|z|^\beta) \qquad \text{as } z \to 0 \text{ in } W \tag{52}$$

*and*

$$u(z) = \mathcal{O}(e^{\gamma \cdot \Re z}) \qquad \text{as } z \to \infty \text{ in } W. \tag{53}$$

*Then the cumulant function*

$$\kappa(\theta) = \int_0^\infty (e^{-\theta x} - 1) u(x) \, dx \tag{54}$$

*admits an analytic continuation from $\{\theta \in \mathbb{C} : \Re \theta > \gamma\}$ to $\{\theta \in \mathbb{C} : \theta \neq \gamma, |\arg(\theta - \gamma)| < \frac{\pi}{2} + \psi\}$ that goes uniformly to 0 as $\theta \to \infty$ in $\{\theta \in \mathbb{C} : \theta \neq \gamma, |\arg(\theta - c)| \leq \frac{\pi}{2} + \psi\}$, where $c > \gamma$ is arbitrary, but fixed.*

*Furthermore, $p(x;t)$ is, for all $x > 0$, an entire function in $t \in \mathbb{C}$ and we have the power series expansion*

$$p(x;t) = \sum_{n \geq 1} u_n(x) \frac{t^n}{n!}, \tag{55}$$

*where*

$$u_n(x) = \frac{1}{2\pi i} \int_{\mathcal{C}} \kappa(\theta)^n e^{\theta x} \, d\theta, \tag{56}$$

*with $\mathcal{C}$ the contour $|\arg(\theta - c)| = \psi$, $\theta = 0$ being passed on the left.*

**Proof.** Let $v(z) = u(z) - az^{-1-\alpha}$ and $\lambda(\theta) = \int_0^\infty e^{-\theta x} v(x) \, dx$. Using

$$\int_0^\infty (e^{-\theta x} - 1)x^{-1-\alpha} \, dx = \Gamma(-\alpha) \cdot \theta^\alpha, \tag{57}$$

we have

$$\kappa(\theta) = a\Gamma(-\alpha) \cdot \theta^\alpha + \lambda(\theta) - \lambda(0), \tag{58}$$

valid for $\Re\theta \geq 0$.

If we let $v_-(x) = v(e^{-i\psi}x)e^{-i\psi}$, then we have the growth estimates $v_-(x) = \mathcal{O}(|x|^\beta)$ as $x \to 0$ and $v_-(x) = \mathcal{O}(\exp((\gamma \cos\psi)x))$ as $x \to \infty$. Thus, the Laplace transform $\lambda_-(\theta) = \int_0^\infty e^{-\theta x} v_-(x) \, dx$ is absolutely convergent for $\Re\theta > \gamma \cos\psi$ and $\lambda_-(\theta) \to 0$ uniformly as $\theta \to \infty$ in $\Re\theta \geq c \cos\psi$; see Doetsch (1950), Satz 4, page 142 and Satz 7, page 171.

Next, we show that $\lambda(\theta) = \lambda_-(\theta e^{-i\psi})$ for real $\theta > \gamma$. Suppose $n \geq 1$ and let us integrate $e^{-\theta z} v(z)$ over the closed contour consisting of a straight line from $n^{-1}$ to $n$, a circular arc from $n$ to $ne^{-i\psi}$, a straight line from $ne^{-i\psi}$ to $n^{-1}e^{-i\psi}$ and a circular arc from $n^{-1}e^{-i\psi}$ to $n^{-1}$ (see Figure 1). By Cauchy's theorem, this integral is zero. The estimates (52) and (53) show that the contributions from the circular arcs vanish as $n \to \infty$ and we obtain

$$\lambda(\theta) = \int_0^\infty e^{-\theta x} v(x) \, dx = \int_0^{e^{-i\psi} \cdot \infty} e^{-\theta z} v(z) \, dz \tag{59}$$

$$= \int_0^\infty e^{-\theta e^{-i\psi} x} v(e^{-i\psi}x) e^{-i\psi} \, dx = \lambda_-(\theta e^{-i\psi}). \tag{60}$$

A similar argument shows that the function $v_+(x) = v(e^{i\psi}x)e^{i\psi}$ has Laplace transform $\lambda_+(\theta) = \int_0^\infty e^{-\theta x} v_+(x) \, dx$, which is absolutely convergent for $\Re\theta > \gamma \cos\psi$ and satisfies $\lambda_+(\theta) \to 0$ uniformly as $\theta \to \infty$ in $\Re\theta \geq c \cos\psi$, and $\lambda(\theta) = \lambda_+(\theta e^{i\psi})$ for real $\theta > \gamma$.

Looking at (58) reveals that the Laplace transform of $p(x;t)$, namely $e^{\kappa(\theta)t}$, is integrable on the vertical line $\Re\theta = c$. Thus, we can use the inversion integral

$$p(x;t) = \frac{1}{2\pi i} \int_{c-i\infty}^{c+i\infty} e^{\kappa(\theta)t + \theta x} \, d\theta. \tag{61}$$

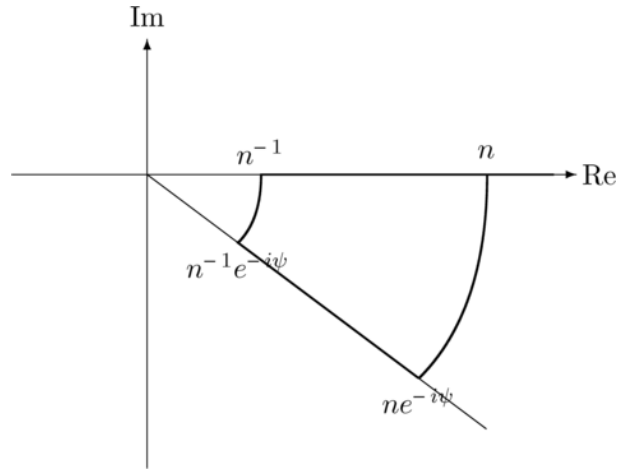

**Figure 1.** Integration contour for the analytic continuation of $\kappa(\theta)$.

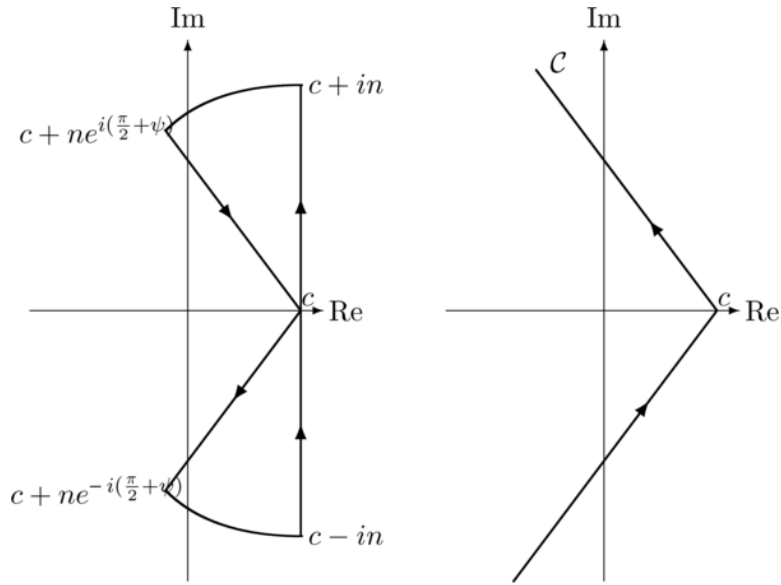

**Figure 2.** The contour used to derive (62) and the final contour $\mathcal{C}$.

Let $n \geq 1$ and consider the integrand $e^{\kappa(\theta)t + \theta x}$ on the closed contour consisting of the vertical line connecting $c - \mathrm{i} \cdot n$ and $c + \mathrm{i} \cdot n$, the circular arc with center $c$ and radius $n$ going from $c + \mathrm{i} \cdot n$ to $c + e^{\mathrm{i}(\pi/2 + \psi)} \cdot n$, the straight lines connecting $c + e^{\mathrm{i}(\pi/2 + \psi)} \cdot n$, $c$ and $c + \mathrm{e}^{-\mathrm{i}(\pi/2+\psi)} \cdot n$ and, finally, the circular arc from $c + \mathrm{e}^{-\mathrm{i}(\pi/2+\psi)} \cdot n$ to $c - \mathrm{i} \cdot n$, see (Figure 2). By Cauchy's theorem, the integral is zero. Again, looking at (58) and the properties of the analytical continuation $\lambda(\theta)$ reveals that the integrand vanishes uniformly on the circular arcs as $n \to \infty$ and, by Jordan's lemma, we conclude that

$$p(x; t) = \frac{1}{2\pi\mathrm{i}} \int_{\mathcal{C}} \mathrm{e}^{\kappa(\theta)t + \theta x} \, \mathrm{d}\theta. \tag{62}$$

On $\mathcal{C}$, the linear term $\theta x$ dominates $\kappa(\theta)t$ as $\theta \to \infty$. Consequently, (62) makes sense for any $t \in \mathbb{C}$, in contrast to (61), where $t > 0$ is required for convergence. We observe that taking $t = 0$ yields

$$\frac{1}{2\pi\mathrm{i}} \int_{\mathcal{C}} \mathrm{e}^{\theta x} \, \mathrm{d}\theta = 0; \tag{63}$$

thus, $p(x; 0) = 0$, according to our convention above. We can differentiate (62) under the integral. Let us consider $h \in \mathbb{C}$ with $|h| \leq 1$ and $t \in \mathbb{C}$ arbitrary. Using (62), we can write the complex difference quotient

$$\frac{p(x; t+h) - p(x; t)}{h} = \frac{1}{2\pi\mathrm{i}} \int_{\mathcal{C}} \frac{\mathrm{e}^{\kappa(\theta)h} - 1}{h} \cdot \mathrm{e}^{\kappa(\theta)t + \theta x} \, \mathrm{d}\theta. \tag{64}$$

Again invoking the asymptotic behavior of $\kappa(\theta)$ as $\theta \to \infty$ on $\mathcal{C}$, we can, by dominated convergence, prove the existence of the complex derivative $\partial p(x; t)/\partial t$ for all $t \in \mathbb{C}$ and the formula

$$\frac{\partial}{\partial t} p(x; t) = \frac{1}{2\pi\mathrm{i}} \int_{\mathcal{C}} \kappa(\theta) \cdot \mathrm{e}^{\kappa(\theta)t + \theta x} \, \mathrm{d}\theta. \tag{65}$$

It follows that $p(x; t)$ is an entire function in $t$ and

$$\frac{\partial^n}{\partial t^n} p(x; t) = \frac{1}{2\pi\mathrm{i}} \int_{\mathcal{C}} \kappa(\theta)^n \cdot \mathrm{e}^{\kappa(\theta)t + \theta x} \, \mathrm{d}\theta. \tag{66}$$

$\square$

# 4. Examples

## 4.1. The positive $\alpha$-stable law

The results for the positive $\alpha$-stable law with $0 < \alpha < 1$ have already been discussed in an informal way in Section 2.1. In view of the simple form of the Lévy density (8), it is straightforward to check that the assumptions for Theorems 3, 6 and 8 are satisfied. Consequently, the results in Section 2.1 can be derived rigorously by our three methods.

## 4.2. The gamma distribution

Suppose $X$ is the gamma process, for which $X_1$ has the law $\Gamma(\nu, \alpha)$ with parameters $\nu = 1$ and $\alpha = 1$. The probability and Lévy densities are

$$p(x; t) = \frac{1}{\Gamma(t)} x^{t-1} \mathrm{e}^{-x}, \qquad u(x) = x^{-1} \mathrm{e}^{-x}. \tag{67}$$

To illustrate our results, we choose the approximation $u_\varepsilon(x) = x^\varepsilon u(x)$. We note that $u_\varepsilon(x) = \Gamma(\varepsilon) p(x, \varepsilon)$ and thus

$$u_{n\varepsilon}(x) = \Gamma(1 + \varepsilon)^n \cdot \varepsilon^{-n} \sum_{k=1}^{n} (-1)^{n-k} \binom{n}{k} p(x; k\varepsilon). \tag{68}$$

So, we are basically in the situation discussed in the third remark of Section 2.5: the convergence of $u_{n\varepsilon}(x)$ to $u_n(x)$ is equivalent to the convergence of the $n$th order difference quotient of (67) at $t = 0$ to the $n$th derivative from the right. The gamma probability density is, for any $x > 0$, an entire function in $t$ and the coefficients in the series expansion (4) are given by

$$u_n(x) = x^{-1} \mathrm{e}^{-x} \sum_{k=0}^{n-1} \binom{n}{k} k! c_k \ln^{n-k-1} x.$$

The numbers $c_k$ arise in the expansion $\Gamma(1 + z)^{-1} = \sum_{n \geq 0} c_n z^n$. They can be expressed explicitly as

$$c_n = \frac{1}{(n-1)!} Y_{n-1}(\gamma, -\zeta(2), 2\zeta(3), \ldots, (-1)^{n-2}(n-2)! \zeta(n-1))$$

with $Y_n$ the complete exponential Bell polynomials as given in Comtet 1970, III.3.c, $\gamma$ the Euler–Mascheroni constant and $\zeta$ the Riemann zeta function.

Let us now imagine that we did not know (67). The function $xu(x)$ is obviously completely monotone and Theorem 4, and so Theorem 3, applies. Let us illustrate the calculation of $u_2(x)$ by formula (33). The tail integral is $U^+(x) = E_1(x)$, where $E_1(x)$ denotes the *exponential integral*; see Abramowitz and Stegun 1992, 5.1.1, page 228. A direct calculation of $(U^+)^{*2}(x)$ is not very explicit. Let us write

$$U^+(x) = V(x) - L(x), \qquad L(x) = \ln x, \qquad V(x) = \ln x + E_1(x). \tag{69}$$

This decomposition is useful because $L(x)$ is simple, while $V(x)$ and its derivatives are integrable at zero. We have $V(0) = -\gamma$ and $V'(0) = 1$. Equation (69) implies

$$(U^+)^{*2}(x) = L^{*2}(x) - 2(L * V)(x) + V^{*2}(x). \tag{70}$$

Next, we observe $L^{*2}(x) = (\ln^2 x - 2\ln x + 2 - \pi^2/6)x$ and thus $[L^{*2}(x)]'' = (2\ln x)/x$. To compute the second derivatives of $(L * V)(x)$ and $V^{*2}(x)$, we can interchange differentiation and convolution by the usual formulas; see Doetsch 1950, 2.14.5, page 115ff.

Namely, we use

$$[V^{*2}(x)]'' = (V')^{*2}(x) + 2V(0)V'(x) \tag{71}$$

and

$$[(V * L)(x)]'' = (V'' * L)(x) + V'(0)L(x) + V(0)L'(x). \tag{72}$$

The convolution integrals on the right-hand side of (71) and (72) can be computed in terms of the exponential integral. Combining the three contributions, those terms cancel and we obtain

$$u_2(x) = 2x^{-1}e^{-x}(\ln x + \gamma), \tag{73}$$

in agreement with (68) above.

Finally, what can we say about Theorem 8 in this case? In its present form, it does not apply since (52) is not satisfied, though formula (56) is correct. The cumulant function is $\kappa(\theta) = -\ln(1 + \theta)$ and

$$u_n(x) = \frac{(-1)^n}{2\pi i} \int_{\mathcal{C}} \ln(1 + \theta)^n e^{\theta x} \, d\theta. \tag{74}$$

Agreement of this formula with (68) can be established by referring to the Hankel contour integrals for the derivatives of $\Gamma(z)^{-1}$ at $z = 1$.

## 4.3. The inverse Gaussian distribution

The inverse Gaussian distribution $IG(\delta, \gamma)$ with $\delta = 1$ and $\gamma = 1$ has a Lévy density of the form

$$u(x) = \frac{1}{\sqrt{2\pi}} x^{-3/2} e^{-x/2} \tag{75}$$

and the probability density is

$$p(x; t) = \frac{1}{\sqrt{2\pi}} t e^t x^{-3/2} e^{-(t^2 x^{-1} + x)/2}. \tag{76}$$

Using the generating function for the Hermite polynomials $H_n(x)$, namely,

$$e^{2xt - t^2} = \sum_{n \geq 0} H_n(x) \frac{t^n}{n!}, \tag{77}$$

we find

$$u_n(x) = \frac{n}{\sqrt{\pi}} 2^{-n/2} x^{-1-n/2} e^{-x/2} H_{n-1}\left(\sqrt{\frac{x}{2}}\right). \tag{78}$$

Let us choose the approximation $u_\varepsilon(x) = \mathrm{e}^{-\varepsilon^2/(2x)}u(x)$. We recognize that this is a multiple of $p(x, \varepsilon)$ and, again, showing the convergence of $u_{n\varepsilon}(x)$ to $u_n(x)$ essentially reduces to a study of the $n$th order difference quotient of $p(x; t)$ at $t = 0$.

Let us look at the second approach, based on the tail integral. Again, $xu(x)$ is completely monotone. We have

$$U^+(x) = \sqrt{\frac{2}{\pi x}} \mathrm{e}^{-x/2} - \mathrm{erfc}\left(\sqrt{\frac{x}{2}}\right), \tag{79}$$

where $\mathrm{erfc}(x)$ is the *complementary error function*; see Abramowitz and Stegun (1992), 7.1.2, page 297. Let us illustrate the computation of $u_3(x)$. By looking at Laplace transforms, we establish $(U^+)^{*3}(x) = 2\sqrt{\frac{2x}{\pi}}\mathrm{e}^{-x/2}(2+x) - (2x^2+6x)\,\mathrm{erfc}(\sqrt{\frac{x}{2}})$ and differentiating $-(U^+)^{*3}(x)$ three times, we obtain

$$u_3(x) = \frac{3}{\sqrt{2\pi}}x^{-5/2}(x-1)\mathrm{e}^{-x/2}, \tag{80}$$

in agreement with (78) above. Finally, Theorem 8 applies, the cumulant function is $\kappa(\theta) = 1 - \sqrt{1+2\theta}$ and we get

$$u_n(x) = \frac{(-1)^n}{2\pi\mathrm{i}} \int_{\mathcal{C}} (1 - \sqrt{1+2\theta})^n \mathrm{e}^{\theta x}\,\mathrm{d}\theta. \tag{81}$$

Agreement of this formula with (78) can be established as follows. First, we substitute $\theta \mapsto (\theta - 1)/2$ and, expanding the integrand by the binomial theorem, we obtain a sum of Hankel integrals of the form (12) with $\alpha = 1/2$, producing a sum of powers of $x$. Using the well-known explicit form of the coefficients of the Hermite polynomials shows (78).

## 4.4. Examples on $\mathbb{R}$: The Meixner and the normal inverse Gaussian distribution

The following example is not covered by the standing assumptions in this paper, as the Meixner distribution (see Schoutens (2003), Section 5.3.10, page 62), is an infinitely divisible distribution on $\mathbb{R}$, not on $\mathbb{R}_{>0}$.

Let us consider the Meixner distribution with parameters $\mu = 0$, $\delta = 1$, $\alpha = 1$ and $\beta = 0$. It has the density

$$p(x; t) = \frac{1}{\pi} 2^{2t-1} \frac{\Gamma(t+\mathrm{i}x)\Gamma(t-\mathrm{i}x)}{\Gamma(2t)}.$$

This expression can be expanded in a series,

$$p(x; t) = \sum_{n \geq 1} u_n(x)\frac{t^n}{n!}, \qquad |t| < |x|,$$

with

$$u_n(x) = \frac{n}{x\sinh(\pi x)} Y_{n-1}(a_1(x), \ldots, a_{n-1}(x)),$$

where

$$a_1(x) = \psi(\mathrm{i}x) + \psi(-\mathrm{i}x) + 2\ln 2 + 2\gamma$$

and

$$a_n(x) = \psi^{(n)}(\mathrm{i}x) + \psi^{(n)}(-\mathrm{i}x) - (-1)^n 2^n (n-1)!\zeta(n) \qquad (n \geq 2).$$

Here, $Y$ again denotes complete exponential Bell polynomials, $\gamma$ is the Euler–Mascheroni constant, $\psi$ is the digamma function and $\zeta$ the Riemann zeta function. Note, however, that here, $p(x;t)$ is not an entire function in $t$ due to the poles of the Gamma function. Thus, we must expect qualitative differences to the cases studied in the present paper.

The normal inverse Gaussian Lévy process with parameters $\mu = 0$, $\delta = 1$, $\alpha = 1$ and $\beta = 0$ has the probability density

$$p(x;t) = \frac{t\mathrm{e}^t}{\pi} \frac{K_1(\sqrt{t^2+x^2})}{\sqrt{t^2+x^2}} \qquad (x \in \mathbb{R}). \tag{82}$$

This expression admits a power series expansion in $t = 0$ for $|t| < |x|$ and the coefficients $u_n(x)$ can be expressed in terms of polynomials and the Bessel K-function. We do not give details here, but note that this result is similar to the Meixner case.

As both the normal inverse Gaussian and the Meixner Lévy processes can be expressed as subordinated Brownian motions, it might be interesting to investigate the power series expansion of a subordinator and the corresponding subordinated process in general. This issue is left to future research.

## 4.5. A bivariate example: The inverse Gaussian–normal inverse Gaussian law

Again, this example is not covered by the assumptions of the rest of this paper as it deals with a bivariate distribution on $\mathbb{R}_{>0} \times \mathbb{R}$. We consider the probability densities

$$p(x,y;t) = \frac{t}{2\pi} \mathrm{e}^t x^{-2} \exp\left[-\frac{1}{2}\left(\frac{t^2+y^2}{x} + x\right)\right] \tag{83}$$

on $\mathbb{R}_{>0} \times \mathbb{R}$. They correspond to an inverse Gaussian–normal inverse Gaussian, or IG-NIG, Lévy process. For properties of this type of law and its origin in a first-passage time problem for a bivariate Brownian motion, see Barndorff-Nielsen and Blæsild (1983), Example 4.1 and Barndorff-Nielsen and Shephard (2001), Example 4.3.

The associated Lévy measure is

$$u(x,y) = \frac{1}{2\pi} x^{-2} \exp\left[-\frac{1}{2}\left(\frac{y^2}{x} + x\right)\right] \tag{84}$$

and the Laplace cumulant function is $\kappa(\theta, \eta) = 1 - \sqrt{1 + 2\theta - \eta^2}$. We can set $u_\varepsilon(x, y) = u(x, y)e^{-\varepsilon^2/(2x)}$ and proceed, just as in Section 4.3 for the IG distribution, to obtain, with (a bivariate extension of) our first method,

$$u_n(x, y) = \lim_{\varepsilon \to 0} \sum_{k=1}^{n} \binom{n}{k} (-1)^{n-k} \varepsilon^{n-k} u_\varepsilon^{*k}(x, y). \tag{85}$$

The resulting series expansion

$$p(x, y; t) = \sum_{n \geq 1} \frac{t^n}{n!} u_n(x, y), \qquad u_n(x, y) = n 2^{-(n-1)/2} x^{-(n-1)/2} u(x, y) H_{n-1}\left(\sqrt{\frac{x}{2}}\right) \tag{86}$$

agrees with the series obtained directly from the explicit expression for $p(x, y; t)$ and the generating function for the Hermite polynomials; see (77).

# Appendix: Auxiliary results

This section provides technical estimates used in the proofs of Theorem 3 and Theorem 6. Several statements look well known and standard. Yet, for rigorous proofs, a careful checking of integrability and differentiability conditions is necessary and, in a few places, delicate, in the present setting.

## A.1. Auxiliary estimates for the cumulant function

The structure of this subsection is as follows. In Lemma A.2, we show that the integrability assumptions (30) of Theorem 3 imply a certain asymptotic behavior of the derivatives of $u(x)$ as $x \to 0$ and as $x \to \infty$. This is used in Lemma A.3 to derive estimates of the derivatives of the cumulant function $\kappa(\theta)$ as $\theta \to \infty$, by partial integration. We then introduce

$$\beta_{nk}(\theta) = B_{nk}(\kappa'(\theta), \ldots, \kappa^{(n-k+1)}(\theta)) \tag{87}$$

and

$$\beta_{nk}^*(\theta) = B_{nk}(|\kappa'(\theta)|, \ldots, |\kappa^{(n-k+1)}(\theta)|), \tag{88}$$

where $B_{nk}$ denotes the partial Bell polynomials, as defined and discussed in Comtet (1970), Section 3, page 144ff. The estimates for the derivatives of $\kappa(\theta)$ are plugged into the Bell polynomials $\beta_{nk}(\theta)$ in Lemma A.4. Using the latter, we obtain estimates for the derivatives of $\lambda_n(\theta)$ as $\theta \to \infty$ in Lemma A.5. Next, Lemma A.6 shows that $u_\varepsilon(x)$ satisfies the assumptions (30) uniformly for $0 < \varepsilon \leq 1$ and thus, applying Lemmas A.2–A.5 to $u_\varepsilon(x)$, gives uniform estimates for the derivatives of $\lambda_{n\varepsilon}(\theta)$ as $\theta \to \infty$. Finally, Lemma A.7 provides a refined estimate for $\beta_{nk}(\theta)$ from the slightly stronger assumptions of Theorem 6.

We are considering distributions on $\mathbb{R}_{>0}$, thus we have $|e^{\kappa(\theta)}| \leq 1$ for $\Re(\theta) \geq 0$. Moreover, $\kappa(\theta)$ is analytic for $\Re(\theta) > 0$ and

$$\kappa^{(n)}(\theta) = (-1)^n \int_0^\infty e^{-\theta x} x^n u(x)\, dx, \qquad \Re(\theta) > 0, n \geq 1. \tag{89}$$

**Definition A.1.** *Suppose $n \in \mathbb{N}$ and $c > 0$. We then say that assumption $\mathcal{A}_n(c)$ holds for $u(x)$ if $u(x)$ is $n$-times continuously differentiable and*

$$\int_0^\infty e^{-cx} x^{k+1} |u^{(k)}(x)|\, dx < \infty, \qquad k = 0, 1, \ldots, n. \tag{90}$$

A consequence of assumption $\mathcal{A}_n(c)$ is, that the Laplace transform $\int_0^\infty e^{-\theta x} x^m u^{(n)}(x)\, dx$ exists for any $m \geq n+1$ and $\Re(\theta) > c$.

**Lemma A.2.** *Suppose $n \geq 2$ and $c > 0$. If assumption $\mathcal{A}_n(c)$ holds, then*

$$\lim_{x \to 0} e^{-\theta x} x^n u^{(n-2)}(x) = 0, \qquad \lim_{x \to \infty} e^{-\theta x} x^n u^{(n-2)}(x) = 0 \tag{91}$$

*for $\Re(\theta) > c$.*

**Proof.** Let $0 < a < b$. Partial integration gives

$$\int_a^b e^{-\theta x} (x^n u^{(n-2)}(x))'\, dx = e^{-\theta x} x^n u^{(n-2)}(x)\big|_a^b + \theta \int_a^b e^{-\theta x} x^n u^{(n-2)}(x)\, dx. \tag{92}$$

We have $(x^n u^{(n-2)}(x))' = n x^{n-1} u^{(n-2)}(x) + x^n u^{(n-1)}(x)$ and see from the integrability assumptions $\mathcal{A}_n(c)$ and $\Re(\theta) > c$ that both integrals in (92) converge to a finite value as $a \to 0$ and $b \to \infty$ (separately). Thus, the limits $\lim_{x \to 0} e^{-\theta x} x^n u^{(n-2)}(x) = \alpha$ and $\lim_{x \to \infty} e^{-\theta x} x^n u^{(n-2)}(x) = \omega$ exist with finite $\alpha$ and $\omega$. But, $\alpha \neq 0$ or $\omega \neq 0$ would imply that $e^{-\theta x} x^{n-1} u^{(n-2)}(x)$ is asymptoically equivalent to $\alpha/x$ as $x \to 0$, respectively to $\omega/x$ as $x \to \infty$. Both properties would contradict the integrability of $e^{-\theta x} x^{n-1} u^{(n-2)}(x)$ that follows, again, from assumption $\mathcal{A}_n(c)$. Thus, we must have $\alpha = 0$ and $\omega = 0$. $\qquad\square$

The following lemma is essentially a reformulation of the well-known fact that an $n$-times differentiable function $f(x)$ with $f^{(k)}(x)$ integrable for $0 \leq k \leq n$ has a Fourier transform $\hat{f}(y)$ which satisfies $\hat{f}(y) = \mathcal{O}(|y|^{-n})$ as $|y| \to \infty$. As we will need uniform growth estimates later, we provide a more detailed statement with explicit bounds.

**Lemma A.3.** *Suppose $n \geq 0$ and $c > 0$. If assumption $\mathcal{A}_n(c)$ holds and we let*

$$L_k(c) = \int_0^\infty e^{-cx} x^{k+1} |u^{(k)}(x)|\, dx, \qquad k = 0, 1, \ldots, n, \tag{93}$$

*then*

$$|\kappa^{(n)}(\theta)| \leq \frac{M_n(c)}{|\theta|^{n-1}}, \qquad \Re(\theta) > c, \tag{94}$$

*where*

$$M_0(c) = L_0(c), \qquad M_n(c) = \sum_{k=0}^{n-1} \binom{n-1}{k} (n)_{n-1-k} L_k(c) \qquad (n \geq 1). \tag{95}$$

**Proof.** For $n = 1$, we have $\kappa'(\theta) = -\int_0^\infty e^{-\theta x} x u(x)\, dx$ and the assertion of the lemma is obvious, namely, $|\kappa'(\theta)| \leq L_0$. For $n = 0$, we can write $\kappa(\theta) = \int_0^\theta \kappa'(\zeta)\, d\zeta$ and the assertion follows, namely $|\kappa(\theta)| \leq L_0 |\theta|$. For $n \geq 2$, we recall $\kappa^{(n)}(\theta) = (-1)^n \int_0^\infty e^{-\theta x} x^n u(x)\, dx$. Let $0 < a < b$. Repeated partial integration gives

$$\int_a^b e^{-\theta x} x^n u(x)\, dx = -\sum_{k=1}^{n-1} \frac{1}{\theta^k} e^{-\theta x} (x^n u(x))^{(k-1)} \Big|_a^b + \frac{1}{\theta^{n-1}} \int_a^b e^{-\theta x} (x^n u(x))^{(n-1)}\, dx \tag{96}$$

*and, by the Leibniz rule, we obtain*

$$(x^n u(x))^{(k-1)} = \sum_{\ell=0}^{k-1} \binom{k-1}{\ell} (n)_{k-1-\ell} x^{n-1-k} x^{\ell+2} u^{(\ell)}(x). \tag{97}$$

From assumption $\mathcal{A}_n(c)$ and Lemma A.2, we conclude, letting $a \to 0$ and $b \to \infty$, that

$$\kappa^{(n)}(\theta) = \frac{(-1)^n}{\theta^{n-1}} \int_0^\infty e^{-\theta x} (x^n u(x))^{(n-1)}\, dx. \tag{98}$$

Using (97), this time with $k = n$, we get

$$(x^n u(x))^{(n-1)} = \sum_{\ell=0}^{n-1} \binom{n-1}{\ell} (n)_{n-1-\ell} x^{\ell+1} u^{(\ell)}(x). \tag{99}$$

This shows that the integral in (98) is bounded by $M_n(c)$. □

**Lemma A.4.** *Suppose $n \geq 0$ and $c > 0$. If assumption $\mathcal{A}_n(c)$ holds, then*

$$\beta_{nk}^*(\theta) \leq \frac{M_{nk}}{|\theta|^{n-k}}, \qquad k = 1, \ldots, n, \tag{100}$$

*where $M_{nk} = B_{nk}(M_1, \ldots, M_{n-k+1})$, the constants $M_1, \ldots, M_n$ are as in Lemma A.3 above, and $B_{nk}$ denote the partial Bell polynomials.*

**Proof.** The Bell polynomials $B_{nk}$ have non-negative coefficients and are therefore increasing functions of each argument. Using the bounds from Lemma A.3, we have

$$\beta_{nk}^*(\theta) = B_{nk}(|\kappa'(\theta)|, \ldots, |\kappa^{(n-k+1)}|) \leq B_{nk}\left( M_1, \frac{M_2}{|\theta|}, \ldots, \frac{M_{n-k+1}}{|\theta|^{n-k}} \right). \tag{101}$$

Using the homogeneity property of the Bell polynomials (see Comtet (1970), Theorem III.3.A), in particular, the last part of the proof, we obtain the desired result.  □

**Lemma A.5.** *Suppose $m \geq 1$, $n \geq 1$ and assumption $\mathcal{A}_n(c)$ holds for $u(x)$ with some $c > 0$, and let $\lambda_n(\theta) = \kappa(\theta)^n$. Then*

$$|\lambda_n^{(m)}(\theta)| \leq \frac{E_{mn}}{|\theta|^{m-n}} \qquad \textit{with } E_{mn} = \sum_{j=1}^{m \wedge n} (n)_j L_0^{n-j} M_{mj}. \tag{102}$$

**Proof.** From the explicit form of Faà di Bruno's formula (see, e.g., Gradshteyn and Ryzhik (2000), (0.43) and Comtet (1970) (Theorems III.3.A and III.4.A)), we get $\lambda_n^{(m)}(\theta) = \sum_{j=1}^{m \wedge n} (n)_j \kappa(\theta)^{n-j} \times \beta_{mj}(\theta)$ and, in conjunction with the estimates from Lemma A.4 we obtain the result.  □

**Lemma A.6.** *Suppose $n \geq 0$, assumption $\mathcal{A}_n(c)$ holds for $u(x)$ with some $c > 0$ and we let*

$$L_k(c) = \int_0^\infty e^{-cx} x^{k+1} |u^{(k)}(x)| \, \mathrm{d}x, \qquad k = 0, 1, \ldots, n. \tag{103}$$

*If we set $u_\varepsilon(x) = e^{-\varepsilon/x} u(x)$ for $x > 0$ and $\varepsilon > 0$, then assumption $\mathcal{A}_n(c)$ holds for $u_\varepsilon(x)$ and we have, for any $\varepsilon > 0$, the uniform bound*

$$\int_0^\infty e^{-cx} x^{n+1} |u_\varepsilon^{(n)}(x)| \, \mathrm{d}x \leq \bar{L}_n(c), \tag{104}$$

*where*

$$\bar{L}_n(c) = \sum_{k=0}^n \sum_{\ell=1}^k \binom{n}{k} \binom{k-1}{\ell-1} \frac{k!}{\ell!} \ell^\ell e^{-\ell} L_{n-k}(c). \tag{105}$$

**Proof.** This follows from

$$[e^{-\varepsilon/x}]^{(k)} = \sum_{\ell=1}^k (-1)^{\ell+k} \binom{k-1}{\ell-1} \frac{k!}{\ell!} x^{-k} \left(\frac{\varepsilon}{x}\right)^\ell e^{-\varepsilon/x} \tag{106}$$

and the inequality $0 \leq x^\ell e^{-x} \leq \ell^\ell e^{-\ell}$ for $x \geq 0$.  □

**Lemma A.7.** *Suppose $m \in \mathbb{N}$, $\alpha \in (0,1)$, $c \in (0,\infty)$ and $u(x)$ is the Lévy density of an infinite activity subordinator which is $m$-times differentiable and such that the functions*

$$v_\ell(x) = \mathrm{e}^{-cx} x^{\ell+1} u^{(\ell)}(x) \qquad (\ell = 0, \ldots, m) \tag{107}$$

*are integrable and their integral modulus of continuity satisfies*

$$\omega^{(1)}(\delta; v_\ell) = \mathcal{O}(\delta^{1-\alpha}) \qquad (\delta \to 0). \tag{108}$$

*We then have, for $n = 0, \ldots, m$,*

$$\kappa^{(n)}(\theta) = \mathcal{O}(|\theta|^{\alpha-n}), \qquad (\Im(\theta) \to \pm\infty) \tag{109}$$

*and, for $\ell = 1, \ldots, m$,*

$$\beta_{m\ell}(\theta) = \mathcal{O}(|\theta|^{\ell\alpha-m}) \qquad (\Im(\theta) \to \pm\infty). \tag{110}$$

**Proof.** From the proof of Lemma A.3 above, we know that

$$\kappa^{(n)}(\theta) = \frac{(-1)^n}{\theta^{n-1}} \sum_{\ell=0}^{n-1} \binom{n-1}{\ell} (n)_{n-1-\ell} \cdot \int_0^\infty \mathrm{e}^{-\theta x} x^{\ell+1} u^{(\ell)}(x) \, \mathrm{d}x \tag{111}$$

for $n = 1, \ldots, m$. Using the assumptions (108) and the well-known relation between the asymptotic behavior of the integral modulus of continuity at zero and the asymptotic growth of the Fourier transform at infinity yields (109) for $n = 1, \ldots, m$. To be more specific, we apply Kawata (1972), Theorem 2.7.4 to the functions $f(x) = \mathrm{e}^{-\Re(\theta)x} x^{\ell+1} u^{(\ell)}(x) I_{(x>0)}$. Note that there is a misprint (not relevant here) in the reference: $f(t)$ should be $\hat{f}(t)$.

The case $n = 0$ follows immediately by using the estimate for $\kappa'(\theta)$ in $\kappa(\theta) = \int_0^\theta \kappa'(\zeta) \, \mathrm{d}\zeta$. Plugging these estimates into (87) and looking at the explicit formula for the Bell polynomials given in Comtet (1970), Theorem III.3.A shows that

$$\beta_{m\ell}(\theta) = \mathcal{O}\Big(\sum |\theta|^{a_1(\alpha-1)+\cdots+a_m(\alpha-m)}\Big) = \mathcal{O}(|\theta|^{\ell\alpha-m}). \tag{112}$$

$\square$

**Lemma A.8.** *Suppose the assumptions for Lemma A.7 hold and*

$$\lambda_n^{(m)}(\theta;t) = \frac{\partial^{m+n}}{\partial\theta^m \partial t^n} \mathrm{e}^{\kappa(\theta)t}. \tag{113}$$

*Then*

$$\lambda_n^{(m)}(\theta;t) = \mathcal{O}(|\theta|^{(m+n)\alpha-m}) \qquad (\Im(\theta) \to \pm\infty). \tag{114}$$

**Proof.** First, we have

$$\lambda_n^{(0)}(\theta;t) = \kappa(\theta)^n \mathrm{e}^{\kappa(\theta)t}. \tag{115}$$

Using (89) and (109) with $n = 0$ shows the claim for $m = 0$. Next, by differentiating (115) $m \geq 1$ times according to Faa di Bruno's formula, we get

$$\lambda_n^{(m)}(\theta;t) = \sum_{\ell=1}^{m} \sum_{j=0}^{\ell} \binom{\ell}{j} (n)_j \kappa(\theta)^{n-j} t^{\ell-j} \mathrm{e}^{\kappa(\theta)t} \beta_{m\ell}(\theta) \kappa(\theta)^n \mathrm{e}^{\kappa(\theta)t}. \tag{116}$$

Using (89), (109) with $n = 0$ and (110), we obtain (114). $\qquad\square$

## A.2. Convolutions and Laplace transforms

In this subsection, we provide further auxiliary results for the proof of Theorem 3. For notational convenience, let us define $V(x) = U^+(x)$ and $V_n(x) = V^{*n}(x)$. First it is shown in Lemma A.9 that the convolution powers $V^{*n}(x)$ exist. We then show that the $V^{*n}(x)$ are $n$-times differentiable and we provide an integral representation in Lemma A.10.

**Lemma A.9.** *Let*

$$V(x) = \int_x^\infty u(y)\,\mathrm{d}y. \tag{117}$$

*Then*

$$V_n(x) = V^{*n}(x), \qquad x > 0, \tag{118}$$

*is well defined for $n \geq 1$ and we have the Laplace transforms*

$$\int_0^\infty \mathrm{e}^{-\theta x} V_n(x)\,\mathrm{d}x = (-1)^n \frac{\kappa(\theta)^n}{\theta^n}, \qquad \Re(\theta) > 0. \tag{119}$$

**Proof.** Let $r > 0$ be arbitrary. A standard argument using the Fubini–Tonelli theorem shows that $\tilde{V}(x) = \mathrm{e}^{-rx} V(x)$ is integrable. Thus, the convolution powers $\tilde{V}^{*n}(x)$ exist for almost all $x > 0$ and are integrable on $\mathbb{R}_{>0}$. As we have

$$\mathrm{e}^{rx} \tilde{V}^{*2}(x) = \int_0^x V(y) V(x-y)\,\mathrm{d}y, \tag{120}$$

this shows that $V(y) V(x-y)$ is integrable on $(0, x)$ and $V^{*2}(x)$ exists for almost all $x > 0$. Repeating the argument shows that the higher convolution powers exist for *almost all* $x > 0$. $\qquad\square$

In Barndorff-Nielsen and Hubalek (2006), a second proof is given, using the theory of convolutions of functions of the class $\mathcal{J}_0$ (from Doetsch (1950)) and it is shown that the convolution powers $V_n(x) = V^{*n}(x)$ actually exist for *all* $x > 0$.

Next, we address the differentiability of $V_n(x) = V^{*n}(x)$. For technical reasons, we work with the derivatives of $x^m V_n(x)$ instead, where $m$ is sufficiently large.

**Lemma A.10.** *Suppose $n \in \mathbb{N}$ and assumption $\mathcal{A}_m(c)$ holds for some $m \geq n + 2$ and $c > 0$. Then $V_n(x)$ is $n$-times differentiable and*

$$V_n^{(k)}(x) = \frac{(-1)^{m+n}}{2\pi \mathrm{i} x^m} \int_{c-\mathrm{i}\infty}^{c+\mathrm{i}\infty} \left( \frac{\kappa(\theta)^n}{\theta^{n-k}} \right)^{(m)} \mathrm{e}^{\theta x} \, \mathrm{d}\theta, \qquad x > 0, \ k = 0, 1, \ldots, n. \tag{121}$$

**Proof.** The result follows basically from Erdélyi *et al.* (1954), IV.4.1 (13), page 130, which implies, that

$$\int_0^\infty \mathrm{e}^{-\theta x} x^m V_n^{(k)}(x) \, \mathrm{d}x = \left( \frac{\kappa(\theta)^n}{\theta^{n-k}} \right)^{(m)}. \tag{122}$$

For a rigorous proof of (121), we have to show that (i) $V_n(x)$ is actually $n$-times differentiable, (ii) that the integrability conditions that allow the application of the quoted rule for the Laplace transform are satisfied and (iii) that the use of the Laplace inversion formula to deduce (121) from (122) is valid. This can be done in an elementary way by induction, with some careful and tedious bookkeeping, incorporating the estimates from Appendix A.1. For details, we refer the reader to Barndorff-Nielsen and Hubalek (2006), Lemmata 23 and 24. □

The result is less obvious than it at first seems. First, the common belief, that convolution increases smoothness is, in general, not true, as the shocking counterexamples in Uludağ (1998) demonstrate. Second, the integrability assumptions to apply the standard theorems on the derivatives of convolutions (such as in Doetsch (1950), I.2.14.5) are typically not satisfied in our setting for $V(x)$, as can be immediately seen from the positive stable example. Third, the application of Erdélyi *et al.* (1954), IV.4.1 (13), page 130 cannot validly be decomposed into an application of Erdélyi *et al.* (1954), IV.4.1 (8), page 129 followed by an application of Erdélyi *et al.* (1954), IV.4.1 (6), page 129.

# Acknowledgements

We thank Ken-Iti Sato, Angelo E. Koudou and an anonymous referee for helpful comments on an earlier version of this paper.